\documentclass[letter,11pt]{article}

\usepackage{comment}
\usepackage[linktocpage=true]{hyperref}
\usepackage[margin=1.2 in, top=1 in, bottom= 1.2 in]{geometry}

\usepackage[english]{babel}
\usepackage{amsfonts}
\usepackage{mathrsfs}
\usepackage[mathscr]{euscript}
\usepackage{bbm}
\usepackage{latexsym}
\usepackage{math dots}
\usepackage{amssymb}
\usepackage{mathtools}

\usepackage{enumitem}
\setlist{nolistsep}
\usepackage{amsthm}
\usepackage{relsize}
\usepackage{nicefrac}
\usepackage[capitalize]{cleveref}

\newtheoremstyle{plain}{3mm}{3mm}{\slshape}{}{\bfseries}{.}{.5em}{}
\newtheoremstyle{definition}{2mm}{2mm}{}{}{\bfseries}{.}{.5em}{}
\theoremstyle{plain}
	
\newtheorem{theorem}{Theorem}

\theoremstyle{definition}

\theoremstyle{plain}
\newtheorem*{namedthm}{\namedthmname}
\newcounter{namedthm}
	

\newcommand{\N}{\mathbb{N}}
\newcommand{\Z}{\mathbb{Z}}
\newcommand{\R}{\mathbb{R}}
\newcommand{\C}{\mathscr{C}}

\newcommand{\eps}{\epsilon}
\newcommand{\mdim}{\overline{D}}

\title{Marstrand type slicing statements in $\Z^{2}\subset \mathbb{R}^{2}$ are false for the counting dimension.}
\author{Aritro Pathak} 
\date{}

\begin{document}
\maketitle

\begin{abstract}
    We show that for $1$ separated subsets of $\R^{2}$, the natural Marstrand type slicing statements are false with the counting dimension that was used earlier by Moreira and Lima and variants of which were introduced earlier in different contexts. We construct a $1$ separated subset $E$ of the plane which has counting dimension $1$, while for a positive Lebesgue measure parameter set of tubes of width $1$, the intersection of the tube with the set $E$ has counting dimension $1$. This is in contrast to the behavior of such sets with the mass dimension where the slicing theorems hold true.
\end{abstract}

\maketitle

\section{Introduction and statement of results.}

\ \ \ \ Fractal properties of subsets of the integer grid in $\R^{d}$ have been studied earlier and notions of dimensions of such subsets have been introduced in different contexts by Fisher \cite{Fisher}, Bedford and Fisher\cite{Bedford}, Lima and Moreira \cite{Moreira}, Naudts \cite{Naudts1, Naudts2}, Furstenberg \cite{Furstenberg}, Barlow and Taylor \cite{Barlow1,Barlow2}, and Iosevich, Rudnev and Uriarte-Tuero\cite{Iosevich}, Glasscock \cite{DG}. The analogies have been drawn from the continuous theory of dimensions, for which Chapter 4 of \cite{Mattila} or Chapter 1 of \cite{Bishop} are standard references.

The mass and counting dimensions of any $1$-separated set $E\in \R^2$ are respectively defined as: 

\begin{equation}
    \mdim(E)=\limsup\limits_{l\to \infty}\frac{\log|E\cap [-l,l]^{2}|}{\log(2l)}, \ D(E)=\limsup\limits_{||C||\to \infty}\frac{\log|E\cap C|}{\log||C||}.
\end{equation}
Here, $C\subset \R^{2}$ is any arbitrary cube with sides parallel to the axes, and $||C||$ denotes the length of each side of $C$.

In \cite{Moreira}, the counting dimension was used in $\Z\subset \R$ to study the growth of certain subsets of $\Z$ with zero upper Banach density. A natural Marstrand type projection theorem is proved there, with the counting dimension in this discrete setting resembling the Hausdorff dimension in the statement of the classical Marstrand projection theorem; see Theorem 1.2 in \cite{Moreira}. Later this was extended by Glasscock \cite{DG} who used the more general notion of the mass and counting dimension, in $\Z^{d}\subset \R^{d}$, and proved analogous projection theorems with the mass dimension as well. 

The natural dual slicing statement with the mass dimension was recently shown to be true by the author \cite{Aritro}. When dealing with the slicing question with a $1$ separated set in $\R^2$, it is natural to work with a width 1 tube $t_{u,v}$ which is explicitly described as:

\begin{equation}
    t_{u,v} = \left\{(x,y) \in \R^2 \ \middle| \ -\frac 1u x + v \sqrt{ 1 + \frac 1{u^2}} < y \leq -\frac 1u x + (v+1) \sqrt{ 1 + \frac 1{u^2}} \right\}.
\end{equation}

This is a tube of width 1. The line perpendicular to this tube and passing through the origin, has slope $u$. Henceforth we call this line $t_{\perp}$. The displacement of the point of intersection of the right edge of $t_{u,v}$ with $t_{\perp}$, is given by the coordinate $v$ .  
This is shown in the figure.

\begin{figure}[h]
\centering
\includegraphics[width=0.4\textwidth]{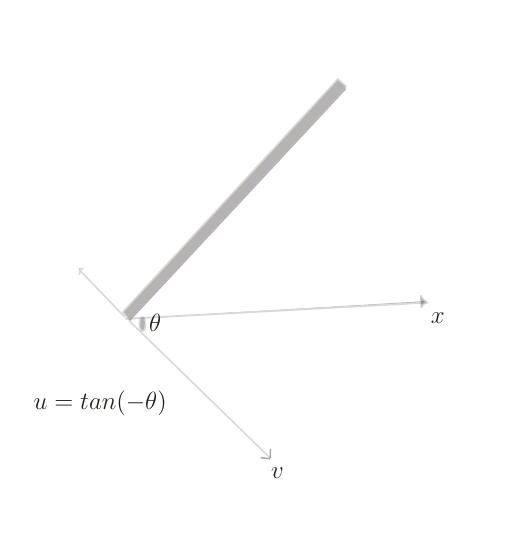}
\caption{The tube $t_{u,0}$, with $v=0$ is shown, and the perpendicular line $t_{\perp}$. A tube $t_{u,v}$ with $v\neq 0$, is a parallel translation of $t_{u,0}$, along the line $t_{\perp}$. }
\label{fig:inksone3}
\end{figure}


In \cite{Aritro}, for the mass dimension in our setting of $1$ separated subsets in $\R^{2}$, with a Tchebysheff and Fubini type argument the slicing statement was first shown to be true in an asymptotic sense, and then it was also shown to be true for Lebesgue almost every slice. One then specializes to sets $A,B \subset \N$ and considers the dimension of the intersection of the broken line $\{(x,y): y=\lfloor \tilde{u}x +\tilde{v} \rfloor, \tilde{u}>0 \}$ with the Cartesian product set $A \times B$. Such a broken line is a tube with vertical cross section of length 1, and thus the width of the tube is less than 1 and so the result follows for such tubes, from the earlier statement for tubes of width $1$.

We state the main slicing result with the mass dimension that was obtained earlier in \cite{Aritro}.

\begin{theorem}\label{thm:strongerslicing3}
Let $E \subseteq  \R^2$ be a $1$ separated set. Then for all $v \in \R$, for Lebesgue-a.e. $u \in \R_+$, 
\[\mdim(E \cap t_{u,v}) \leq \text{max} (0, \mdim(E)-1).\]
\end{theorem}

As a corollary to this, we obtained the main slicing result stated below, upon integrating over the $v$ coordinate.

\begin{theorem}\label{thm:slicing}
Let $E \subseteq  \R^2$ be a $1$ separated set. Then in the Lebesgue sense, for almost every tube $t_{u,v}$ of width $1$, slope $u$, and displacement $v$ along the projecting line, we have that $\mdim(E \cap t_{u,v}) \leq \text{max} (0, \mdim(E)-1)$.
\end{theorem}

In this paper we show that the corresponding results with the counting dimension are false. This is what one would expect, since the counting dimension of every slice can be high if there are cluster of points in every slice at very sparsely separated locations, with the number of points in each cluster growing to infinity. This can happen even if the actual set $E$ has low counting dimension.

We state our results for the counting dimension below:

\begin{theorem}\label{thm:strongerslicing1}
\begin{enumerate}
\item For any $\epsilon>0, u_{0}\in \R$, there exists $E \subseteq  \Z^2$ , such that for all $v \in (-\epsilon, \epsilon)$,  
\[\overline{D}(E \cap t_{u_0,v})=D(E \cap t_{u_0,v})= D(E)=\overline{D}(E)=1\]
\item For any $\epsilon>0, u_{0}\in \R$, there exists $E \subseteq  \Z^2$ , such that for all $v \in (-\epsilon, \epsilon)$,  
\[D(E \cap t_{u_0,v})= D(E)=1, \overline{D}(E \cap t_{u_0,v})=\overline{D}(E)=0. \]
\end{enumerate}
\end{theorem}

The construction in the part (1) of this theorem is actually trivial; an analogous version  \cref{thm:strongerslicing1} is also true for the mass dimension, cf. Example 1 in Section 4 of \cite{Aritro}. This just consists of the set of the form $\{(x,y\in \Z^{2}: -\eps/2 < x\leq \eps/2), y\geq 0\}$.

The more interesting case of \cref{thm:strongerslicing2} is to construct the set $E$ that is meagre in terms of the mass dimension, in part (2) of the above theorem .

Next we state the main theorem of this paper, 
\begin{theorem}\label{thm:strongerslicing2}
For any $\epsilon>0, v_{0}\in \R$, there exists $E \subseteq  \Z^2$, such that for all $u \in (-\epsilon, \epsilon)$,  
\[D(E \cap t_{u,v_0})= D(E)=1.\]
\end{theorem}

This shows the distinction from the behavior of the mass dimension, where for any one-separated subset $E\subset\R^2$, we have for any $v\in \R$, for Lebesgue-a.e. $u\in \R$, $\overline{D}(E\cap t_{u,v})\leq \text{max}(\overline{D}(E)-1,0)$, which is the content of Theorem 9 in \cite{Aritro}.

Without loss of generality, in proving \cref{thm:strongerslicing2}, it would be enough to construct the set $E$ in the case of $v_0 =0$. The case of $v_0 \neq 0$ can be dealt with, with minor adjustments. 

The same construction of \cref{thm:strongerslicing2}, will provide the counterexample to the statement of the Marstrand slicing theorem with the counting dimension. We state this separately as a theorem.

\begin{theorem}\label{thm:strongslicing}
For any $\epsilon>0$, there exists $E \subseteq  \Z^2$ , such that for all $(u,v) \in (-\epsilon, \epsilon)\times (-\epsilon,\epsilon)$,  
\[D(E \cap t_{u,v})= D(E)=1.\]
\end{theorem}

In other words, we can construct a set $E\subset \Z^{2}$ such that the slices of $E$ by the tubes parametrized by a positive Lebesgue measure set have exceptionally high counting dimension.

In fact this construction of \cref{thm:strongslicing} for any fixed $\epsilon$, immediately gives us the following result: 

\begin{theorem}\label{thm:weakreal}
\emph{There exists a set $E \subseteq \N^2$, so that the set $U$ of parameters $(u,v)$, with $u,v\in \R$, so that}
\[D \big(E \cap  t_{(u,v)} \big) = D(E)=1 .\] \emph{is such that $\lim\limits_{M\to \infty} \frac{|U \cap [-M,M]^{2}|}{(2M)^{2}}>0$. }
\end{theorem}
One compares this with the behavior of the mass dimension, Theorem 6 in \cite{Aritro}.

In the next section, as an illustration we begin with a standard example that shows that the slicing result is sharp with the mass dimension, where we have a set so that every slice in a cone has mass dimension $\frac{1}{2}$ while the set so constructed has mass dimension $\frac{3}{2}$. This is a set where every slice in the cone has counting dimension 1, while the set itself has a subsequence of grid of points where the number of points in each grid grows in a two dimension sense, to infinity, and hence is itself of counting dimension 2. After that we construct a set where we have a growth of points in a linear 'diagonal' manner within the cone which reduces the counting dimension of the set to 1, while every slice in the cone still has counting dimension 1. 

We also show that the proof of \cref{thm:weakreal} essentially follows from the construction used in \cref{thm:strongslicing}. This is the strongest slicing statement we can make with the counting dimension. 

We remark that the analysis here again remains the same if in place of $\Z^{2} \subset \R^{2}$, we considered the grid $\delta\Z^{2}$ with separation $\delta$, and considered tubes of width $\delta$. 

\textbf{Notation:} Throughout the rest of the paper, we use the term `cluster' of points to mean a subset of the form $T_{a,b}=\{(x,y)\in \Z^{2}: a_1\leq x\leq a_2, b_1\leq y\leq b_2 \}$ for some  integers $a_1\leq a_2, b_1\leq b_2$. We would also talk about clusters of points within a cone $\C$, or a cluster of points within a tube $t_{(u,v)}$ in which case it will be understood to mean $T_{a,b}\cap \C$, or $T_{a,b}\cap t_{(u,v)}$ respectively. When talking about such clusters, we will not enumerate the set of values of $a_1,a_2,b_1,b_2$ as above, and it should be understood from the context that we mean a cluster of this type. Also, in several places, we use the term 'level' to indicate the group of clusters of the set $E$  within the cone $C$, which are at approximately the same height from the origin (such as shown in Figure 2).

\section{Slicing results with the counting dimension.}

 Consider a cone of arbitrary small angular width $\theta$, centered around the $y-$axis, with vertex at the origin and pointing up. Let us call this cone $\mathcal{C}$. Let $l_1,l_2$ be the two limiting lines of the cone. In this case, for some large $ k_0 >0$, \footnote{In order to construct a set with mass dimension 3/2, we chose the numbers $2^{2^{k}}$ instead of $2^k$ since if we have a range of $k$ values from 1 to some $k_0$, the levels are so sparse that only the last level $k_{0}$ is relevant when counting the points for the mass dimension till height $2^{2^{k_0}}$. With the scaling $2^{k_0}$ we would need to add up all the points in all the lower levels as well.} for each $k\geq k_0$, we fill the annular region inside the cone between the heights $2^{2^{k+1}}$ and $2^{2^{k+1}}+ 2^{2^{k}}$ with all the points belonging to $\mathbb{Z}^{2}$ within the annular region; i.e. we consider the intersection of $\mathcal{C}$ with the region $\{(x,y)\in \Z^2: 2^{2^{k+1}}\leq y\leq 2^{2^{k+1}}+ 2^{2^{k}}   \}$. This example was already considered in \cite{Aritro}. This is a set of mass dimension $3/2$ with each of the tubes within $\C$ having mass dimension $1/2$. However this is also a set that has counting dimension 2: since for all $k\geq k_0$ above the height $2^{2^{k+1}}$ we effectively have a square of dimension $2^{2^{k}} \times 2^{2^{k}}$ that is filled with points of the integer grid. Considering this sequence of cubes, as $k$ is taken to infinity, this implies that we have a set of counting dimension 2. Moreover, each of the tubes has an intersection with $\sim 2^{2^{k}}$ many points just above the height $2^{2^{k+1}}$ and so as $k$ is taken to infinity, each of the tubes also have counting dimension 1. 

Retain the setting of the cone of the previous paragraph. In order to construct a set $E$ supported within this cone $\C$, with counting dimension $D(E)>1$ such that we have a family of tubes $\{t_{(u,0)}\}$ for $u$ belonging in some small interval around the origin, with each tube in the family having an exceptionally high counting dimension greater than $(D(E)-1)$, we have to ensure that the set of points of $E$ do not accumulate in two dimensional clusters of growing length. In the examples below, we show how to construct a set $E$ that grows linearly. We can't locate within this set any growing two dimensional sub-sequence that is growing to infinity. This is done while ensuring that  each of the ray of tubes has a growing sub-sequence of points that ensures that it still has counting dimension 1. That would prove \cref{thm:strongerslicing2}.
 
\bigskip 

Now we construct the set that proves \cref{thm:strongerslicing1}. 


\begin{proof} [Proof of \cref{thm:strongerslicing1}, (2)]
 Let $u_{0}=0$, and consider a set $V$ of values of $v$, with $\mu(V)=wk$, where $w>1$ is a positive number and $k$ is any arbitrary positive integer and $\mu$ is the Lebesgue measure \footnote{We can clearly make this work for any $\eps$ so that the statement of the theorem is satisfied}. In this case, the set $E$ is contained within the semi-infinite strip of width $wk$, $\{(x,y)\in \Z^2: - \frac{wk}{2} <x\leq \frac{wk}{2}, y\geq 0 \}$.
 
 The idea here is similar to the one used prior to this theorem, and we start with a cluster of the integer grid of width $w$, height $n_{1}$, placed at the bottom left corner of this semi-infinite tube. We place the next cluster of the same width and height on the top right corner of the previous block, and do this all the way till we put the last block adjacent to the right edge of this tube, at the height $H_{1}:=kn_{1}$.
 
 Now we consider the height $h_{2}:=e^{H_{1}}$ and repeat the same process at this level with the clusters placed diagonally as before from the bottom left to the top right corner, with now the blocks of height $n_{2}:=n_1 +1$ and width $w$. Thus at the end we reach the height $H_{2}:=h_{2}+wn_{2}$. Inductively we repeat the process where at each step $H_{m}:=h_{m}+wn_{m}$, where $h_{m}:=e^{H_{m-1}}$ and $n_{m}:=n_{m-1}+1$. Like before, it is clear that this is a set where between the heights $h_{m}$ and $H_{m}$ we have for the purpose of the counting dimension, in effect a straight block of width $w$ and length $kn_{m}\to \infty$ as $m\to \infty$. It is clear this set $E$ has counting dimension exactly 1, when taking boxes of lengths $kn_{m}\to \infty$ as $m\to \infty$ at the $m$'th level of the construction of the set $E$. Moreover, every single tube with $v$ parameter within $V$ has between $n_{m}$ and $2n_{m}$ many points at the $m$'th level, and thus also clearly a tube of counting dimension exactly equal to 1.
 
In this case, it is clear that  $H_{m}> \underbrace{e^{e^{\iddots }}}_{m \ \text{times}}$ and till that height we have $wk(n_1 + (n_1 +1)+(n_1+2)+\dots+(n_1 +m))=\big(wkm(n_1 + \frac{m+1}{2})  \big)$ many points in E.  Thus clearly this set $E$ has mass dimension 0. It's similarly also clear that every width 1 tube with parameter $(u_{0},v)$, $v\in V$, has mass dimension 0,
\end{proof}
\bigskip

\begin{figure}[h]
\centering
\includegraphics[width=0.5\textwidth]{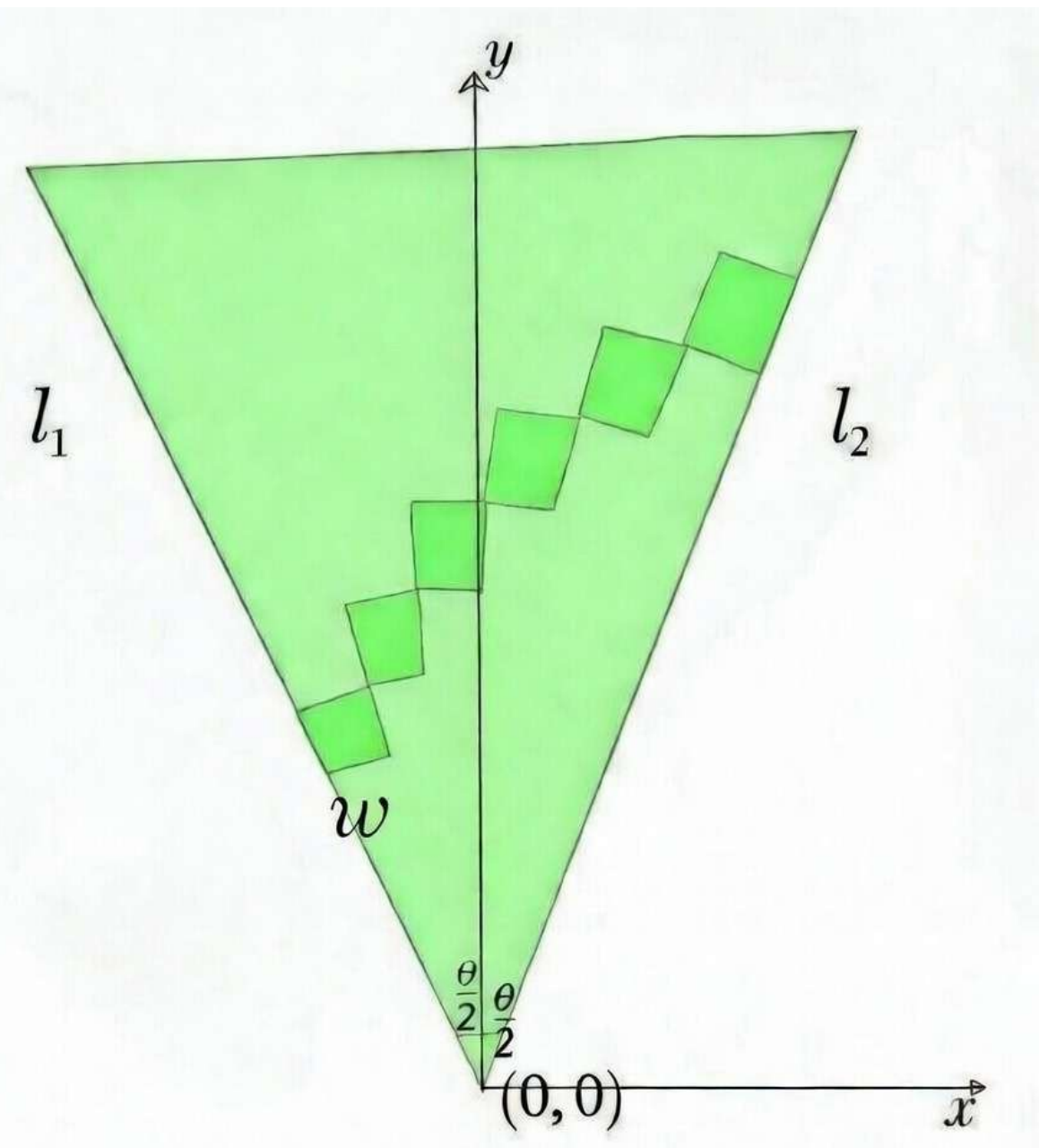}
\caption{The cone showing the clusters at any specific level of the construction of the set $E$, growing 'diagonally' from one edge of the cone to the other. Every level of $E$ consists of such a diagonal sequence of points, and only one level of $E$ is shown in the figure. The width of each cluster is the same $w$, while the height of each cluster at the level $k$, is given by the sequence $n_k \to\infty$ as $k\to \infty$, which for simplicity can be chosen to be $n_k =k, \ \forall k\in \Z$.}
\label{fig:inksone3}
\end{figure}

Now we prove \cref{thm:strongerslicing2}.

\begin{proof}[Proof of \cref{thm:strongerslicing2}]
Without loss of generality, consider a cone with the vertex at the origin, and pointing up, symmetric about the y axis, and with total angle $\theta$. Call the left edge of the cone $l_1$, and the right edge of the cone $l_2$. Thus, the angles within the cone are parametrized by the coordinate $u\in (-\theta/2, \theta/2)$. Here, all the tubes within the cone have the coordinate $v=0$, and we have an interval of $u$ values of width $\theta$. 

We fix an integer $w\geq1$, an initial length $n_{1}\geq 1$. Adjacent to $l_1$, we put a cluster  (or more generally any maximal $1$ separated set) of width $w$ and length $n_{1}$,\footnote{Since $n_1\geq 1, w\geq 1$ this block would contain at least one point.} just above the arc at height $h_{1}$. This subtends the angle $(-\frac{\theta}{2},-\frac{\theta}{2}+\frac{w}{h_{1}})$ at the origin, of width $\frac{w}{h_1}$.
At the next step, we put a block of width $w$ just above the arc at height $h_1 +n_1$ so this subtends the angle $(-\frac{\theta}{2}+\frac{w}{h_{1}},-\frac{\theta}{2}+ \frac{w}{h_1}+ \frac{w}{h_1+n_1})$ of width $\frac{w}{h_1+n_1}$. At the $k$'th step, we put a block of width $w$ at the height $h_1 +(k-1)n_1$ so that it subtends the angle $(-\frac{\theta}{2}+\sum\limits_{i=0}^{k-2} \frac{w}{h_{1}+i\cdot n_{1}},-\frac{\theta}{2}+\sum\limits_{i=0}^{k-1}\frac{w}{h_{1}+i\cdot n_1})$ of width $\frac{w}{h_1 +(k-1)\cdot n_1}$. Clearly the growth of the total angle is like that of the harmonic series, and eventually after some step $K_1$ we reach the other end of the cone at the angle $\theta/2$. In this process, it is clear that the width of the last cluster adjacent to $l_2$, at the height $H_{1}:=h_1 +(K_1-1)n_1$ could be less than $w$.

Consider the next height $h_2:=e^{H_{1}}$, and $n_{2}:=n_1 +1$ while the width $w$ remains fixed. Now the growth of points here is in a manner similar to the growth in the first level, and we begin the growth again from the left edge $l_1$ of the cone starting at the height $h_{2}$. At the $k$'the stage in this level, we put a block of width $w$ at the height $h_2 +(k-1)n_{2}$ which subtends the angle $(-\frac{\theta}{2}+\sum\limits_{i=0}^{k-2}\frac{w}{h_2 +i\cdot n_{2}}, -\frac{\theta}{2}+\sum\limits_{i=0}^{k-1}\frac{w}{h_{2}+i\cdot n_{2}})$ of width $\frac{w}{h_{2}+(k-1)\cdot n_{2}}$. Again, after a finite number $K_2$ of steps, we would hit the right edge $l_2$ of the cone, where again at the last end, the width of the integer grid just above the height $H_2:=h_2 +(K_{2}-1)n_2$ is less than or equal to $w$.

For each $m\geq 3$, we would inductively define $h_{m}:=e^{H_{m-1}}$ and $n_{m}:=n_{m-1}+1$ \big(Thus $n_{m}=n_{1}+(m-1)$\big) and thus have a growth of points beginning at the left edge $l_1$ at height $h_m$ and continuing on to height $H_{m}$ at the last $K_{m}$'th step of the iteration in this $m$'th level, where this last cluster of the integer grid is put adjacent to the right edge $l_2$ of the tube.

This is a set where every tube parametrized by $(u,0)$ with parameter $u\in (-\frac{\theta}{2},\frac{\theta}{2})$, has counting dimension exactly 1; since for each of these tubes, at the $m$'th level, there is a cluster within the tube where the number of points is between $n_{m}$ and $2n_{m}$. Further, by construction, these clusters are exponentially separated. 

But the set $E$ so constructed is itself also of counting dimension exactly 1, since at each particular $m$'th level, we have a set of clusters, each of width $w(1+o(1))$, each of length $n_1 + (m-1)$, the $k$'th one lying just above and to the right of the $(k-1)$'th cluster. If we had a square of length $K_{m}\cdot (n_{1}+(m-1))$ intersect the set $E$ at this appropriate height, then for the purpose of the counting dimension this is equivalent to having a vertical straight line of length $\approx K_{m}\cdot(n_1 +(m-1))$ within this square. As $m\to \infty$, the lengths of these boxes  go to infinity, and we have a set of counting dimension exactly 1.  
\end{proof}

\bigskip

 Note that the mass dimension of this set considered in \cref{thm:strongerslicing2} is 1, and the mass dimension of all the tubes is exactly 0. To see this, note that the height $H_{m}\sim h_{m}.e^{\frac{m\theta}{w}}$ since, the total angle covered is the fixed $\theta \sim \int_{0}^{K_{m}} \frac{w dx}{h_{m}+x n_{m}}$ and where $H_{m}=h_{m}+ (K_{m}-1)n_{m}$ and thus $K_{m}\sim \frac{h_{m}}{n_{m}}\big(e^{\frac{\theta n_{m}}{w}} -1\big)$ and thus $H_{m}=h_{m}+(K_{m}-1)n_{m}\sim h_{m}e^{\frac{\theta m}{w}}$. Thus the growth of the set is equivalent to having a straight line of fixed width $w$, starting at height $h_{m}$ and ending at approximately the height $h_{m}e^{\frac{m\theta}{w}}$, and thus of length about $h_{m}(e^{\frac{\theta m }{w}}-1)$. By construction, we also have that $h_{m}=e^{H_{m-1}}$ for all $m\geq 2$, and so that this height range is from $e^{H_{m-1}}$ to $e^{(H_{m-1}+\frac{\theta \cdot m}{w})}$ and where $H_{m-1}\sim h_{m-1}e^{\frac{\theta\cdot (m-1)}{w}}>> \frac{\theta\cdot m}{w}$ when $m$ is sufficiently large enough. Thus for the purpose of the mass dimension, we consider the appropriate box of length $e^{(H_{m-1}+\frac{\theta \cdot m}{w})}$, and we have a `line' of width $w$ of points from the height $e^{H_{m-1}}$ to $e^{(H_{m-1}+\frac{\theta\cdot m}{w})}$ within this box, while below the height $H_{m-1}$ there are some $O(w\cdot H_{m-1})$ many points. Thus the mass dimension is given by the limit below, as $m\to \infty$.
 
 \[ \frac{\log(O(w\cdot H_{m-1})+ w(e^{(H_{m-1}+\frac{\theta\cdot m}{w})}-e^{H_{m-1}}))}{\log(e^{H_{m-1}+\frac{\theta\cdot m}{w}})} \to \frac{\log(w)+(H_{m-1}+\frac{\theta\cdot m}{w})-\log(1-e^{-\frac{\theta\cdot m}{w}})}{H_{m-1}+\frac{\theta\cdot m}{w}}\to 1 .\]
 
 Also note that each of the tubes has $w\cdot n_{m}\sim wm$ number of points between the heights $h_{m}$ and $H_{m}$ where $h_{m}=e^{H_{m-1}}>e^{m-1}$ for $m\geq 2$, and thus it is easily shown that  each of the tubes has 0 mass dimension. 

 We notice the similarity in the construction of the sets considered in \cref{thm:strongerslicing1}(2) and \cref{thm:strongerslicing2}, in that the growth of the points is in a diagonal manner in both cases. However, the set $E$ constructed in the proof of Theorem 3(2) has mass dimension $0$ while the set constructed in the proof of Theorem 4 has mass dimension $1$.
 
 \bigskip

Now with this same set $E$, we can further show that the natural Marstrand slicing statement with the counting dimension is false. 

\begin{proof}[Proof of \cref{thm:strongslicing}]
The construction follows from the previous example, by considering the same set $E$, withing the cone $C$. The set of tubes is parametrized by $\{(u,v): -\beta \leq u\leq \beta, -\epsilon \leq v \leq \epsilon \}$. With minor modifications to the argument in the previous theorem, we conclude that the set $E$ is such that the intersection of every tube with parameter in the range $(u,v) \in (-\beta ,\beta)\times (-\epsilon,\epsilon)$, with $E$, has counting dimension exactly 1. The minor change here is that these tubes parameterized by the $u\neq 0$ parameters intersect the clusters of $E$ 'obliquely', while still having counting dimension 1.
\end{proof}

Now the proof of \cref{thm:weakreal} is also immediate with the same set $E$. 

\begin{proof}[Proof of \cref{thm:weakreal}]
Again consider the same set $E$, within the same cone. Consider all the tubes $t_{(u,v)}$ whose right edges are lines that intersect the set $E$ with slopes in the range $( \cot \theta/2,\infty)$, and thus the slope $m(t_{\perp})$ of the line $t_{\perp}$ is so that $0\leq m(t_{\perp})\leq \tan (\theta/2)$ (The same is true of tubes whose right edges have negative slopes in the range $(-\tan (\theta/2),-\infty)$). In this case, any one of these tubes intersects all the `levels' of $E$ parametrized by $k$ large enough. It is easy to see that except for at most the first level it intersects, such a tube intersects both the top and bottom boundary layers of all the other levels, and thus these tubes have counting dimension exactly 1.

The set of such tubes is parametrized by the set $S=\{(u,v)\in \Big( (0,\tan (\theta/2))\times (0,\infty     )\Big) \cup \Big((-\tan (\theta/2),0)\cup (-\infty,0)\Big) \} $ and thus $\lim_{M\to \infty} \frac{ S \cap [-M,M]^{2}}{(2M)^{2}} >0$ and so the weak asymptotic form of the slicing theorem is also false for the counting dimension.
\end{proof}

\section{Further questions.}
\begin{enumerate}
     \item While we show that a positive Lebesgue measure set of tubes of width 1 have exceptionally large dimension, the question of constructing sets such that a set of broken lines $\{(x,y): y=\lfloor ux+v \rfloor \}$ with positive Lebesgue measure in the parameter space $(\tilde{u},\tilde{v})$, all have exceptionally large counting dimension, is not resolved here. The broken lines are sets that have vertical cross section 1, but the actual width is smaller than unity, and so there is no guarantee that a large number of broken lines have exceptionally high dimension, let alone a set of broken lines parametrized by a positive Lebesgue measure set in the $(u,v)$ space. Such a question will likely require a much finer study of the intersection patterns of tubes of arbitrary small width with the integer grid $\Z^{2}\subset \R^{2}$. 
     
     \item While we have constructed our sets in $\R^2$, there should be a natural way to extend these results in higher dimensions, and while the construction in \cref{thm:strongerslicing1} would be extended to higher dimension in an obvious way, we would likely have to be a bit more careful in constructing the cones in dimensions 3 and higher while proving the results analogous to \cref{thm:strongerslicing2}.

     \item We aim to study the question of sets of exceptions to the Marstrand type slicing statement of \cite{Aritro} for the mass dimension, as well as the Marstrand type projection theorems proved in \cite{DG, Moreira} and the dimensions of the parameter set of tubes that violate the inequality in the Lebesgue a.e. sense. This will likely involve ergodic and diophantine approximation techniques, and draw on the literature for the study of exceptional sets for the classical Marstrand projection and slicing theorems.

\end{enumerate}
\section{Acknowledgements:} The author thanks Daniel Glasscock for introducing the author to the slicing problem in the integer lattice grid in $\R^{2}$.

\end{document}